\documentclass[twoside,12pt]{article}
\usepackage{amsmath}
\usepackage{amssymb}
\usepackage{amscd}
\usepackage{amsthm}%+

\numberwithin{equation}{subsection}
 
\theoremstyle{plain}

\newtheorem{Le}[equation]{Lemma}

\theoremstyle{remark}

\theoremstyle{definition}

\newcommand{\R}{\mathbb R}
\def\ga{\gamma}

\def\ra{\rightarrow}

\def\e{\emph}
\def\i{\infty}
\def\p{\partial}
\def\b{\begin}

\begin{document}

\title{    \flushleft{\bf{Growth of relatively hyperbolic groups}}      }
%\author{\flushleft{Xiangdong Xie}}
\date{  }
\maketitle

\vspace{-5mm}

\noindent 
Xiangdong Xie\newline
Department of Mathematical Sciences,
 University  of Cincinnati,\newline
Cincinnati, OH 45221, U.S.A.\newline
Email:   {  xxie@math.uc.edu}

%\noindent
%{\bf{Mathematics Subject Classification(2000).}} Primary 53C23, 51F99, 57M20; 
% Secondary  53C70, 57M60.   
% \newline
%{\bf{Keywords.}} Tits boundary, Tits metric, CAT(0), 2-complex, geodesic.

%\keywords    test  \endkeywords
%\subjclass  51  \endsubjclass

\pagestyle{myheadings}

\markboth{{\upshape Xiangdong Xie}}{{\upshape  Uniform exponential growth }}

%\footnotetext{test}
      
%\keywords    test  \endkeywords
%\subjclass  51  \endsubjclass
\vspace{3mm}

\noindent
{\small {\bf Abstract.}
We show that a  relatively hyperbolic group
 either is virtually cyclic or has  uniform exponential growth.}

%\tableofcontents  

\vspace{3mm}

%20F34 Fundamental groups and their automorphisms 
%20F65 Geometric group theory 
%20F67 Hyperbolic groups and nonpositively curved groups
%20F99 None of the above, but in this section
%20E07 Subgroup theorems; subgroup growth
%53C20 Global Riemannian geometry, including pinching
%53C23 Global topological methods 
%57M20 Two-dimensional complexes
%57M60 Group actions in low dimensions
%(1) (2) (a) (b) (i) (ii)
%20F69 Asymptotic properties of groups

%$\{x_i\}_{i=1}^\i$  converges to $\xi\in \ol{X}$ 
% an index two subgroup

\vspace{3mm}

\setcounter{section}{1}
\setcounter{subsection}{0}

\subsection{Introduction}

Let $G$ be a finitely generated group, and $S$ a finite generating set. 
Denote by $d_S$ the word metric on $G$ with respect to $S$ and 
$\beta(G,S,k)$ the number of elements of $G$  a $d_S$ distance at most 
$k$ from the identity.  
The exponential  growth rate of $G$ with respect to $S$ is 
$\omega(G, S)=\lim_{k\ra \infty}(\beta(G,S,k))^{1/k}$. Notice  that 
the limit exists due to the submultiplicativity: 
$\beta(G,S, m+n)\le \beta(G,S,m)\beta(G,S,n)$.  We say 
$G$  has \e{exponential growth} if $\omega(G, S)>1$
  for  some (hence any) finite  generating set $S$; and we say   $G$ has
    \e{uniform exponential growth}
  if $\omega(G)=\inf_S \omega(G, S)>1$, where $S$ varies over all 
finite generating sets of 
  $G$.

 J. Wilson constructed groups that have exponential growth but not uniform exponential 
 growth (\cite{W}).  On the other hand,  among the following classes of groups 
exponential growth  implies uniform exponential 
 growth:  linear groups over fields with zero characteristic (\cite{EMO}), 
 hyperbolic groups (\cite{K}), one-relator groups (\cite{GD}), 
solvable groups (\cite{O2}),  and geometrically finite groups acting on 
 pinched Hadamard manifolds (\cite{AN}).  In this paper we show  that 
the same is true for relatively hyperbolic groups.

\b{Th}\label{main}
{Let $G$ be a relatively hyperbolic group. Then $G$ either 
  has uniform
 exponential growth or has a finite index cyclic subgroup.
}
\end{Th}

Relatively hyperbolic groups are generalizations of Gromov hyperbolic groups.
 Typical examples of relatively hyperbolic groups include Gromov
 hyperbolic groups,  fundamental groups of finite volume real hyperbolic manifolds 
 and groups acting properly and cocompactly on spaces with isolated flats (\cite{HK}).
There are 5 different but equivalent definitions for relatively hyperbolic groups:
 one due to  each of Gromov (\cite{G}),  B. Farb (\cite{F})  and  D.  Osin (\cite{O1}), 
 and 
two due to B. Bowditch (\cite{B}).  We shall use B. Bowditch's definition of 
relatively hyperbolic groups as geometrically finite groups. 
  The main ingredient in our proof is B. Bowditch's theorem on the existence of
 an invariant   collection of  disjoint horoballs (see Proposition \ref{inva}).

The usual way for proving uniform exponential growth is as follows:
  for any  finite generating set 
 $S$ of $G$, find two elements $g_1, g_2\in G$ with word length bounded independent of $S$,
  such that $<g_1, g_2>$ is free with basis $\{g_1, g_2\}$. 
  We shall use the same  strategy.
Let $G$ be a relatively hyperbolic group
  and $S$ a   finite generating set.  
 %of $G$, $d_S$ denotes the word metric on $G$ with respect to $S$.
For any positive integer $n$, let $S(n)=\{g\in G:  d_S(g, id)\le n\}$.  
 Notice that $S(n)$ is also a finite generating set of $G$, and
 $S(n)\subset S(m)$ if $n\le m$. In particular, $S\subset S(n)$ for all $n\ge 1$.

We prove uniform exponential growth in two steps:\newline
Step 1:  There exists some positive integer $n_0$ with the following property.
For any finite generating set $S$ of $G$,   $S(n_0)$ 
contains a hyperbolic element.  See Proposition \ref{p3}.\newline
Step 2:  There exists some positive integer $k_0$ with the following property.
 If a finite generating set $S$ of $G$ contains a hyperbolic element,
 then there  are $g_1, g_2\in S(k_0)$ such  that 
$<g_1, g_2>$ is free with basis $\{g_1, g_2\}$.
See  Corollary  \ref{vc1} and Proposition \ref{p2}.

 %some $y\in S$, and 
 %some integer $k$, $1\le k\le k_0$, 
  %  such that $g_1:=x^k$ and $g_2:=yg_1y^{-1}$   generate 
 %a free group with basis $\{g_1, g_2\}$.  
 %See Proposition \ref{p1} and Proposition \ref{p2}. 

It should be noted that  uniform exponential growth for 
    geometrically finite groups acting on 
 pinched Hadamard manifolds has been established by  
R. Alperin  and   G. Noskov (\cite{AN}). Some of our arguments are similar to theirs. 
  Theorem \ref{main} was conjectured by C. Drutu in \cite{D}.

{\bf{Notation:}}  
 For a metric space $Y$, any subset $A\subset Y$ and any $r\ge  0$, 
 we denote by $N_r(A)=\{y\in Y:  d(y, A)\le r\}$
the  closed $r$-neighborhood of $A$.  For any $p, q\in Y$, $pq$ denotes a geodesic
  segment 
 between $p$ and $q$, although it is not unique in general.

%\subsection{Hyperbolic element case}\label{fold}

\subsection{Relatively hyperbolic groups as geometrically finite groups} \label {subcom}

 Here we recall the notion of geometrically finite groups and B. Bowditch's result
 on the existence of an invariant system
   of disjoint horoballs (see Proposition \ref{inva}).
 The reader is referred to \cite{B} 
for more details. 

Suppose that $M$ is a compact metrizable  topological space  and  a group 
$G$  acts by homeomorphisms on $M$. We say that $G$ is a \e{convergence group}
   if the induced action on the space of distinct triples is properly discontinuous. 
Let $G$ be a convergence group on  $M$.   A point $\xi\in M$ is a \e{conical limit point}
 if there exists a sequence in $G$, $\{g_n\}_{n=1}^\i$,  and  two 
points $\zeta\not=\eta\in M$, such that $g_n(\xi)\ra \zeta$ and 
$g_n(\xi')\ra \eta$ for all $\xi'\not=\xi$.  
An element $g\in G$ is a \e{hyperbolic element}
    if it has infinite order and fixes exactly two points in $M$.  
  A subgroup  $H<G$ is parabolic if it is infinite, 
  fixes a point $\xi\in M$, and  contains no hyperbolic elements.  
  %The nontrivial elements in a parabolic subgroup are called \e{parabolic elements}.
  In this case, the fixed point of $H$ is unique and is
     referred to  as a parabolic point. 
The nontrivial elements in a parabolic subgroup are called \e{parabolic elements}.
 A parabolic point $\xi$ is a \e{bounded
 parabolic point} if its stabilizer $Stab(\xi):=\{g\in G: g(\xi)=\xi\}$ acts properly and cocompactly on 
$M\backslash \{\xi\}$.     A convergence group $G$ on $M$ is a 
\e{geometrically finite group} if each point of $M$ is either a conical limit point or a bounded parabolic point.

\b{Def}\label{rel}
{A group $G$ is \e{hyperbolic relative}
to  a family of finitely generated subgroups ${\cal G}$, if it acts properly discontinuously by isometries,
 on a proper  geodesic hyperbolic space $X$, such that the induced action on $\p X$  is of convergence,
 geometrically finite,  and such that the maximal parabolic subgroups are exactly 
 the elements of $\cal{G}$.}
\end{Def}

By definition a relatively hyperbolic group is a geometrically finite group. 
 A. Yaman (\cite{Y}) proved that if $G$ is a geometrically finite group on a perfect
 metrizable compact space $M$,  and the maximal parabolic subgroups are finitely generated,
  then $G$ is hyperbolic relative to the  family of maximal parabolic subgroups. 
%By a result of P. Tukia (\cite{T}, Theorem 1B),  

Now let $X$ be a $\delta$-hyperbolic geodesic  metric space
for some $\delta\ge 0$. Recall that  each geodesic triangle
 in $X$ is $\delta$-thin, that is, 
 each edge is contained in the $\delta$-neighborhood of the 
union of the other two edges.
 %We will adopt B. Bowdutch's notation.
%Suppose $x, y\in \R$. We write $x\preceq y$ if $x \le y+ c$, where $c$ is a constant
  %depending on $\delta$. 
  %and $\xi\p X$.  
 Below we shall denote by $c=c(\delta)$ a constant that depends only on
$\delta$.  
Let $\xi\in \p X$.  
 A (not necessarily continuous) function $h: X\ra \R$ is a \e{horofunction}
  about $\xi$   if  there are constants $c_1=c_1(\delta)$, 
 $c_2=c_2(\delta)$ such that:      
 if  
  $x, a\in X$    and $d(a, x\xi)\le c_1$ for some geodesic 
ray $x\xi$ from $x$ to $\xi$, then $|h(a)-(h(x)+d(x,a))|\le c_2$.
A closed subset $B\subset X$  is a \e{horoball} about $\xi$  
if there is a  horofunction   $h$ about  $\xi$  
  and a constant $c=c(\delta)$
 such that $h(x)\ge -c$ for all $x\in B$,  and $h(x)\le c$ for all $x\in X\backslash B$. 
Note that $\xi$ is uniquely determined by $B$, and we refer to it as the center of the horoball. %Note also that a horoball is quasi-convex.
%Recall that a subset $A\subset Y$ of a metric space is quasi-convex if there is some 
 %$\lambda\ge 0$ such that for any $a_1, a_2\in A$, every 
%geodesic segment $a_1a_2$  lies in $N_\lambda(A)$. 
% and for any $a\ge 0$, the closed regular neighborhood
%$N_a(B)$ is also a horoball. 

\b{Prop} \e{(B. Bowditch, Proposition 6.13 in \cite{B})}\label{inva}
{Let $G$ be a relatively hyperbolic group,  and $X$ a proper  hyperbolic 
 geodesic metric  space  that $G$ acts upon as in Definition \ref{rel}.  Let $\Pi$ be the set of all 
 bounded parabolic points in $\p X$. Then  $\Pi/G$ is finite. Moreover, for any $r\ge 0$,
 there is a  collection of horoballs ${\cal{B}}=\{B_\xi: \xi\in \Pi\}$ indexed by
 $\Pi$ with the following properties:\newline
\e{(1)} $\cal{B}$ is  $r$-separated, that is, $d(B_\xi,B_\eta)\ge r$ for all $\xi\not=\eta\in \Pi$;\newline
\e{(2)}  ${\cal{B}}$ is $G$-invariant, that is, $g(B_\xi)=B_{g(\xi)}$ for all $g\in G$ and
 $\xi\in \Pi$;\newline
\e{(3)} $Y({\cal{B}})/G$ is compact, where 
$Y({\cal{B}})=X\backslash \bigcup_{\xi\in \Pi}\text{int} (B_\xi)$.}
 \end{Prop}

%It follows from  Proposition  \ref{inva} 
%that there are a finite number of pairwise nonconjugate subgroups   $\{H_1, \cdots,  H_k\}$  of  $G$ such that  
% $\{gH_ig^{-1}: 1\le i\le k, g\in G\}$ is  the family of maximal parabolic subgroups, 
%that is, the family of stabilizers at all the bounded parabolic points. 
 %In this case, we also say $G$ is hyperbolic relative to 
%$\{H_1, \cdots,  H_k\}$. 

\subsection{Axes of hyperbolic elements} \label{axis}

 In this section we study how a hyperbolic element
  \lq\lq translates" its  \lq\lq axes". 

Let $G$ and $X$ be as in Definition \ref{rel}. That is, 
 $X$ is a proper $\delta$-hyperbolic geodesic space for some $\delta\ge 0$, 
 $G$  acts properly discontinuously by isometries on $X$, 
% on a proper hyperbolic space $X$,
 such that the induced action on $\p X$  is of convergence  and 
 geometrically finite. %  and such that the maximal parabolic subgroups are exactly 
 %the elements of $\cal{g}$.

Recall  that,  in a $\delta$-hyperbolic space $X$, any two complete geodesics that
have the same endpoints in $\p X$ have Hausdorff distance at most $2\delta$. 
%$\delta$-thin
  For   a hyperbolic element  $\ga\in G$,
  let  $\ga_+$  and  $\ga_-$    be the attracting and repelling   fixed
points of $\ga$ in  $\p X$ respectively. 
We shall call any complete geodesic 
 with $\ga_+, \ga_-$ as endpoints  an axis of $\ga$,  and denote by $A_\ga$ the union of all
 axes of $\ga$.   Note that 
 $\ga$ may have many different axes and an axis of $\ga$
  in general is not invariant under $\ga$. However, for any axis $c$ of $\ga$,
   $\ga(c)$  is  also an axis of $\ga$ and hence 
  the Hausdorff distance between $c$ and $\ga(c)$ is at most $2\delta$.

 Proposition \ref{inva} implies that there is a 
$200\delta$-separated $G$-invariant collection of horoballs $\cal{B}$ centered at the parabolic points  such that $Y({\cal{B}})/G$ is compact. 

%Let
% ${\cal{B}}'=\{N_{49\delta}: B\in {\cal{B}}\}$.  Then 
%${\cal{B}}'$ is also a $2\delta$-separated  $G$-invariant
%collection of horoballs. 

\b{Le}\label{n0}
{There exists a positive integer $k_1$ with the following property:
  for any   infinite order  element 
$\ga\in G$    %any axis $c$ of $\ga$, 
and any $x\in Y({\cal{B}})$,  
%$x\in c\cap Y({\cal{B}})$,  
there is some $k$, $1\le k\le k_1$ such that
  $d(\ga^k(x), x)\ge 200\delta$.}
\end{Le}

\b{proof}  Notice that   the 
action of $G$ on $Y({\cal{B}})$ is  properly discontinuous  and cocompact.
  It follows that  there is some integer 
$k_1\ge 1$ such that for any $x\in Y(\cal{B})$,  
  the cardinality of $\{g\in G: d(x, g(x))\le 200\delta\}$ is less than 
 $k_1$.  In particular, for any  $x\in Y({\cal{B}})$ and any 
infinite order  element
   $\ga\in G$, there is some 
$k$, $1\le k\le k_1$ such that   $d(\ga^k(x), x)\ge 200\delta$.

\end{proof}

\b{Cor}\label{vc1}
{%There exists a positive integer $k_1$ with the following property.
  If a   finite generating set  $S$ of $G$  contains a hyperbolic element,
    then  there is a hyperbolic 
 element $\ga\in S(k_1)$  and some 
$x\in A_\ga$  such that $d(\ga(x), x)\ge 200\delta$, where
 $k_1$ is the constant in Lemma \ref{n0}.}
\end{Cor}

%any axis 
%$c$  of $\ga$, there is some  $x\in c$ and 
%some 
%$k$, $1\le k\le k_0$ such that   $d(\ga^k(x), x)\ge 200\delta$.}
%\end{Cor}

\b{proof} Let $g\in S$ be a hyperbolic element and $c$ an axis of $g$.
 It follows from the definition of a horoball that 
$c$ is not contained in  any horoball. 
Since $\cal{B}$ is a disjoint collection of horoballs and c is connected, 
 we have $c\cap Y(\cal{B})\not=\emptyset$. 
Let $x\in c\cap Y(\cal{B})$. By Lemma \ref{n0}, there is some 
$k$, $1\le k\le k_1$ such that   $d(g^k(x), x)\ge 200\delta$.
 Notice that $g^k\in S(k_1)$ is hyperbolic and $c$ is also an axis of 
$g^k$.  

\end{proof}

For a complete geodesic $c$ in $X$, we define a map $P_c: X\ra c$
  as follows:   for any $x\in X$,  
  let  $P_c(x)\in c$ be a point with $d(x, P_c(x))=d(x,c)$. 
Note that  for any two points $x_1, x_2\in c$ with 
$d(x, x_1)=d(x, x_2)=d(x,c)$, we have $d(x_1, x_2)\le 4 \delta$.

\b{Le}\label{or1}
{Let $c, c': \R\ra X$ be two geodesics with 
$c(+\i)=c'(+\i)$ and $c(-\i)=c'(-\i)$.  
Let $a, b\in \R$ with $b\ge a+8\delta$,  and $a'$, $b'\in \R$ 
   be  determined by
$c'(a')=P_{c'}(c(a))$
  and $c'(b')=P_{c'}(c(b))$.   Then $(b'-a')\ge (b-a)-4\delta$.}
\end{Le}

\b{proof}
Recall   that 
the Hausdorff distance between $c$ and $c'$ is at most $2\delta$. 
It follows from triangle inequality that
$|b'-a'|\ge (b-a)-4\delta$.
  So either 
$(b'-a')\ge (b-a)-4\delta$ or $(b'-a')\le  -(b-a)+4\delta$. 
Suppose $(b'-a')\le  -(b-a)+4\delta$.    
    Let  $x_i=c(b+i\delta)$, $i=0, 1, \cdots$,  and  $y_i=P_{c'}(x_i)$. Then 
 $y_0=c'(b'), y_1, \cdots, $  is a sequence of points
  on $c'$ with $d(y_i, y_{i+1})\le 5\delta$. Since 
   $(b'-a')\le  -(b-a)+4\delta\le -4\delta$  and  $y_i\ra c(+\i)$,  
there is some $i\ge 1$ such that 
  $d(y_i, c'(a'))\le 2.5 \delta$.  Triangle inequality implies that 
$(b+i\delta)-a=d(x_i, c(a))\le d(x_i, y_i)+d(y_i, c'(a'))+d(c'(a'), c(a))\le
  6.5 \delta$, contradicting the fact that $b\ge a+8\delta$.

\end{proof}

\b{Le}\label{order}
{Let $g\in G$ be a hyperbolic element, and $c$ an axis of $g$.
   Suppose  $x\in c$  is a point with  
$d(g(x), x)\ge 20\delta$.  Then for any integers 
$i< j$,  the point $P_c(g^j(x))$ lies between 
 $P_c(g^i(x))$  and $g_+$.}
\end{Le}

\b{proof}
Denote $x_i=P_c(g^i(x))$.  Note that $x_i\ra g_+$ as $i\ra +\i$. 
It suffices to show that for any $i$, $x_i$ lies between 
$x_{i-1}$ and $x_{i+1}$. 
Set $t=d(g(x), x)$.  
Triangle inequality implies that $t-4\delta\le d(x_i, x_{i+1})\le t+4\delta$.
Suppose  $x_i$  does not lie between 
$x_{i-1}$ and $x_{i+1}$  for some $i$.  
Then  $d(x_{i-1}, x_{i+1})\le 8\delta$
  and $d(g^{i-1}(x), g^{i+1}(x))\le 12 \delta$.  Since $g$ is an isometry, 
$d(x, g^2(x))\le 12\delta$. We shall show that $d(x, g^2(x))\ge 2t-16\delta$,
 which is a  contradiction.

Let $\sigma: \R\ra X$ be the parameterization of $c$ with 
 $\sigma(0)=x$ and $\sigma(+\i)=g_+$, $\sigma(-\i)=g_-$.
 Then  $g\circ \sigma$ is the  parameterization of $g(c)$
 with $g\circ \sigma(0)=g(x)$  and  $g\circ \sigma(+\i)=g_+$,
$g\circ \sigma(-\i)=g_-$.
%Note that $c$ and $g(c)$ have the same endpoints and the same direction. 
%We parameterize $c$ and $g(c)$ such that $x=c(0)$ and $g(x)=g(c)(0)$. 
Since $t-4\delta\le d(x_1, x)\le t+4\delta$, we have 
$x_1=\sigma(b)$ with either $-(t+4\delta)\le b\le -(t-4\delta)$ or
$t-4\delta\le b \le t+4\delta$.  We consider the case 
$-(t+4\delta)\le b\le -(t-4\delta)$, the other case  can be handled  similarly. 
Let $y_1=P_c(g(x_1))$. Then $y_1=\sigma(a)$ for some $a\in \R$. 
Now Lemma \ref{or1}  applied to $\sigma$, $g\circ \sigma$ and the projections of
 $g(x)$, $g(x_1)$ onto $c$ 
implies that $b-a\ge 0-b-4\delta\ge t-8\delta$. It follows that 
 $d(y_1, x)=0-a\ge -b+t-8\delta\ge 2t-12\delta$. By triangle inequality
 $d(x, g^2(x))\ge d(x, y_1)-
d(y_1, g(x_1))-d(g(x_1), g^2(x))\ge 2t-12\delta-2\delta-2\delta=2t-16\delta$.
%We are done.

\end{proof}

%For any hyperbolic element  $\ga\in G$,  $A_\ga$ denotes the 
%union of all axes of $\ga$.

\b{Le}\label{pars}
{Let $g\in G$ be a hyperbolic element. 
Suppose there is a point  $x\in A_g$ with $d(x, g(x))\ge 200\delta$. Then 
$|d(y, g(y))-d(z, g(z))|\le  40 \delta$   for all
 $y,z\in A_g$.} 

\end{Le}

\b{proof}
Let $c$ be an axis of $g$ that contains $x$. %Denote $x_i=P_c(g^i(x))$. 
We first show $|d(y, g(y))-d(x, g(x))|\le 16\delta$ for all $y\in c$.
Fix $y\in c$. Denote $t=d(x, g(x))$, 
$x_i=P_c(g^i(x))$,  $x'_i=P_c(g(x_i))$ and $y'=P_c(g(y))$. 
Then $d(x_i, g^i(x))\le 2\delta$  and $d(x_i, x_{i-1})\ge t-4\delta$.
By Lemma \ref{order}  
there is some $i$ such that 
 $y$ lies between $x_i$ and $x_{i+1}$,  and   $x_{i+1}$
    lies between $x_i$ and $g_+$.  % and $x_{i+1}$ lies between $x_i$ and $g_+$. 
Triangle inequality implies 
$|d(y', x'_{i-1})-d(y,x_{i-1})|= 
|d(y', x'_{i-1})-d(g(y),g(x_{i-1}))| \le 4\delta$
 and $d(x_i, x'_{i-1})\le d(x_i, g^i(x))+
d(g^i(x), g(x_{i-1}))+d(g(x_{i-1}),  x'_{i-1}))\le 6\delta$. 
It follows that
 $|d(y', x_i)-d(y,x_{i-1})| \le 10\delta$. 
In particular, $d(y', x_i)\ge d(y,x_{i-1})-10\delta=
d(y, x_i)+d(x_i, x_{i-1})-10\delta\ge d(y, x_i)+t-4\delta-10\delta> d(y, x_i)$.
%It follows that $d(P_c(g(y)), x_i)\ge d(P_c(g(y)), P_c(g(x_{i-1}))-d(P_c(g(x_{i-1})), x_i)
%\ge d(y, x_{i-1})-4\delta-d(P_c(g(x_{i-1})), x_i)\ge d(y, x_{i-1})-10\delta
%=d(y, x_i)+d(x_i, x_{i-1})-10\delta\ge d(y, x_i)+t-4\delta-10\delta> d(y, x_i)$,
 %that is, $d(P_c(g(y)), x_i)> d(y, x_i)$. 
On the other hand, Lemma \ref{or1}
  applied to $g(c)$, $c$ and the projections of $g(y)$, $g(x_{i-1})$ 
  onto  $c$ implies that
 $y' $ lies between $x_i$ and $g_+$. 
It now follows from $d(y', x_i)> d(y, x_i)$ that 
 $y$ lies between $x_i$ and 
$y' $.  Now   $|d(y, y')-d(x_{i-1},x_i)|=
|[d(x_i, y)+d(y, y')]-[d(x_{i-1},x_i)+d(x_i,y)]|
=|d(x_i, y')-d(x_{i-1}, y)|\le 10\delta$.
%|d(x_i, P_c(g(y)))-d(x_{i-1}, y)|
%=|[d(x_i, y)+d(y, P_c(g(y)))]-[d(x_{i-1},x_i)+d(x_i,y)]|\le 
%|d(y, P_c(g(y)))-d(x_{i-1},x_i)$
 Finally $|d(y, g(y))-d(x, g(x))|= |d(y, g(y))-d(g^{i-1}(x), g^i(x))|
\le |d(y, y')-d(x_{i-1},x_i)|
+ d(g(y), y')+d(g^{i-1}(x), x_{i-1})+d(g^i(x), x_i)
\le 16\delta$.

Now fix any $y\in A_g$.  Set $y'=P_c(y)\in c$. 
Then $|d(y, g(y))-d(x, g(x))|\le 
|d(y', g(y'))-d(x, g(x))|+d(y, y')+d(g(y'), g(y))\le 20\delta$.
 By triangle inequality,  
we have 
$|d(y, g(y))-d(z, g(z))|\le 40\delta$  for all $y, z\in A_g$.

\end{proof}

\subsection{Free subgroups} \label{frees}

The goal in this section is to find two short hyperbolic elements that generate 
a  free group:

\b{Prop}\label{p2}
{Let $G$ be a relatively hyperbolic group,  and $X$ a proper  $\delta$-hyperbolic 
 geodesic metric  space  that $G$ acts upon as in Definition \ref{rel}.
Then 
there exists a positive integer $k_2$ 
with the following property.
  If   $S$ is a      finite generating set    of $G$ 
 and 
$s\in S$ is a hyperbolic element such that 
 $d(s(x), x)\ge 200\delta$ for some $x\in A_s$,  then 
 there are $g_1, g_2\in S(k_2)$ such that $<g_1, g_2>$ is 
 free with basis $\{g_1, g_2\}$.}

\end{Prop}

\subsubsection{Ping-Pong lemma} \label{ping}

The proof of Proposition \ref{p2} is based on the following Ping-Pong lemma.
%For any group  $G$ and   any  $g_1, g_2\in G$,  $<g_1, g_2>$  denotes the 
 % subgroup generated by $g_1$  and $g_2$. 
%We recall the following well-known lemma.

\begin{Le}\label{pingpong}
{Let $G$ be a group acting on  a  set $X$, and  $g_1$, $g_2$   two  elements of $G$. 
If $X_1$, $X_2$ are disjoint subsets of $X$ and for all $n\not=0$, $i\not=j$,
  $g_i^n(X_j)\subset X_i$, then the subgroup $<g_1, g_2>$ is free 
  with  basis $\{g_1, g_2\}$.}

\end{Le} 

We will apply the Ping-Pong Lemma in the following setting.

\begin{Le}\label{appping}
{Let $X$ be a metric space and $g_1, g_2$ isometries of $X$.
Let $B_i\subset A_i\subset X$ \e{(}$i=1,2$\e{)} be subsets of $X$,  and 
$p_i: X\ra A_i$  a map. Denote $X_i=p_i^{-1}(A_i-B_i)$. 
 Then the assumptions in the Ping-Pong lemma are satisfied if the following 
  conditions  hold:\newline
\e{(1)} $X_1\cap X_2=\emptyset$;\newline
\e{(2)}  $g_i^n(p_i^{-1}(B_i))\subset X_i$  
for all $n\not=0$, $i=1,2$.}
\end{Le}

\qed

In our case,  $g_i$ is a  hyperbolic element, $A_i$ is  an axis of $g_i$,
 $B_i\subset  A_i$  is a segment, and $p_i=P_{A_i}$.

%We shall show that there exists some constant $N_0$ with the following property:
 % for any finite generating set $S$ of $G$, 
 %there are $x, y\in S$ and some  integer $k$, $1\le k\le N_0$ 
 %such that $g_1:=x^k$ and $g_2:=yg_1y^{-1}$ satisfy 
 %the assumptions in the above lemma. Compare with the remark in the Introduction. 

\subsubsection{Two hyperbolic elements} \label{twohe}

In this section we prove Proposition \ref{p2}. For this,  we 
shall find two short hyperbolic elements and segments in 
 their axes such that the two  conditions in Lemma \ref{appping}
 are satisfied.

Let $\cal{B}$  be a fixed
$200\delta$-separated $G$-invariant collection of horoballs 
   centered at the parabolic points. Recall (see the proof of Lemma \ref{n0}) that 
there is an  integer 
$k_1\ge 1$ such that for any $x\in Y(\cal{B})$,  
  the cardinality of $\{g\in G: d(x, g(x))\le 200\delta\}$ is less than 
 $k_1$.

\b{Le}\label{twoax}
{Let $g\in G$  be  a    hyperbolic element  such that 
 there is some $x\in A_g$  with $d(x, g(x))\ge 200\delta$,   and  $\ga\in G$
 an element such that  $\{\ga(g_+), \ga(g_-)\}\cap \{g_+, g_-\}=\emptyset$. 
%that does not   share  any   fixed points  with $g$.
 Let $c_1$ and $c_2$ be  axes of $g$ and $h:=\ga g\ga^{-1}$
  respectively. 
Let $a=\min\{t\in \R: d(c_1(t), c_2)\le 2\delta\}$
 and $b=\max\{t\in \R: d(c_1(t), c_2)\le 2\delta\}$, and denote $y=c_1(a)$, 
$z=c_1(b)$.  Then  $d(y,z)\le 3k_1 d(x, g(x))$.}
%  where $k_1$  is the integer in Lemma \ref{n0}.}

\end{Le}

\b{proof}  Denote  $t=d(x, g(x))$.  
Suppose $d(y,z)> 3k_1 t$. 
Let $y'=c_2(a')\in c_2$, $z'=c_2(b')\in c_2$ with $d(y,y'), d(z,z')\le 2\delta$. 
 It is easy to see from the $\delta$-thin condition that the Hausdorff distance between 
 $yz$ and $y'z'$ is  at most $4\delta$.    
After  replacing $h$ by its inverse if necessary, we may assume 
$g$ and $h$ translate in the \lq\lq same direction", that is, if both $c_1$ and $c_2$ are parameterized from the repelling fixed point toward the attracting fixed point, then 
 $b'>a' $ (we have $b>a$ by definition).

 Since $h$ is a conjugate of $g$ or its inverse,
 Lemma \ref{pars}  implies that the inequality
$|d(p, g^i(p))-d(q, h^i(q))|\le 40\delta$ holds
for all
  $i\ge 1$,  and all  $p\in A_g$,  $q\in A_h$. 
In particular, 
$|d(y, g^i(y))-d(y', h^i(y'))|\le 40\delta$
for all $i\ge 1$.

For $i\ge 0$, 
  let $y_i=P_{c_1}(g^i(y))$  and  $y'_i=P_{c_2}(h^i(y'))$.
 Lemma \ref{order} implies that 
 $y_{i+1}$ lies between $y_i$ and $c_1(+\i)$.     Since  
$d(y_i, y_{i+1})\le d(y, g(y))+4\delta\le d(x,g(x))+40\delta+4\delta\le 2t$,
  we have   $y_i\in yz$ for all $1\le i\le k_1$.
Similarly $y'_i\in y'z'$ for all $1\le i\le k_1$.
Hence for  $1\le i\le k_1$
there exists $z_i\in yz$ with $d(z_i, y'_i)\le 4\delta$.
 %Note that $z_i$ lies between $y$ and $c_1(+\i)$.

We first assume $y\in Y(\cal{B})$.      
%Then $|d(y, g^i(y))-d(y', h^i(y'))|\le 108\delta$. 
Triangle inequality implies 
%$|d(y', y'_i)-d(y, z_i)|\le d(y,y')+d(z_i, y'_i)\le 8\delta$   and 
$|d(y, y_i)-d(y', y'_i)|\le |d(y, g^i(y))-d(y', h^i(y'))|
 + d(y'_i, h^i(y'))+d(y_i,  g^i(y))\le 40\delta+2\delta+2\delta=44\delta$
 and $|d(y', y'_i)-d(y, z_i)|\le d(y,y')+d(z_i, y'_i)\le 6\delta$.
It follows that $d(y_i, z_i)=|d(y, y_i)-d(y, z_i)|\le 50\delta$ and 
$d( g^i(y), h^i(y'))\le d(g^i(y),y_i)+d(y_i, z_i)+d(z_i, y'_i)+d(y'_i, h^i(y'))
\le 2\delta+ 50\delta+4\delta+2\delta=58\delta$. 
 Now  for any $i$, $1\le i\le k_1$,  we have  $d(y, h^{-i}\circ g^i(y))=
d(h^i(y), g^i(y))\le d(h^i(y), h^i(y'))+d(h^i(y'), g^i(y))
=d(y, y')+d(h^i(y'), g^i(y))\le 60\delta$. 
 Since $y\in Y(\cal{B})$, it follows from 
 the definition of $k_1$  that  
%the cardinality of 
%$\{g\in G: d(y, g(y))\le 200\delta\}$ is less than $k_1$.
  there are $i\not=j$, $1\le i,j\le k_1$  with  
$h^{-i}\circ g^i=h^{-j}\circ g^j$. Consequently $h^{j-i}=g^{j-i}$,
  contradicting the fact that $g$ and $h$ do not share any 
fixed points  in $\p X$.

Now suppose $y\notin Y(\cal{B})$. Since $Y(\cal{B})$ is $G$-invariant,
$g(y)\notin Y(\cal{B})$. There are horoballs $B_1, B_2\in \cal{B}$ with
$y\in B_1$ and $g(y)\in B_2$.  Note $B_1\not= B_2$, otherwise
 $g(B_1)=B_1$ and so $g$ fixes the center of $B_1$, contradicting the fact that $g$ 
 is a hyperbolic element. %Let $y_1=P_{c_1}(g(y))$.
Since $\cal{B}$ is $200\delta$-separated  and $d(y_1, g(y))\le 2\delta$, 
    there is some $p\in c_1$ between $y$ and $y_1$ such that 
$p\in Y(\cal{B})$. Now  we run the argument from the pervious paragraph using 
 $p$ instead of  $y$.

\end{proof}

\b{remark}\label{notwos}
{The proof of Lemma \ref{twoax} actually shows that
for any $\ga\in G$, either $\{\ga(g_+), \ga(g_-)\}=\{g_+,  g_-\}$ or 
$\{\ga(g_+), \ga(g_-)\}\cap \{g_+,  g_-\}=\emptyset$  holds. }

\end{remark}

Now let $S$ be a finite generating set of $G$   and  
 suppose $s\in S$ is a hyperbolic element  such that $d(s(x), x)\ge 200\delta$ for some 
  $x\in A_s$. 
 If $\{s'(s_+), s'(s_-)\}=\{s_+,  s_-\}$ for all $s'\in S$, then 
$A_{s}$  is invariant under $G$. In this case, 
$G$ is $2$-ended and hence is  virtually infinite cyclic.

We now assume there is some $\ga\in S$ with 
$\{\ga(s_+), \ga(s_-)\}\cap \{s_+,  s_-\}=\emptyset$.
Let $c$ be an axis of $s$. 
Set $g_1=s^{10k_1}$ and $g_2=\ga g_1\ga^{-1}$,
 $A_1=c$, $A_2=\ga(c)$, $p_1=P_c$ and $p_2=P_{\ga(c)}$. 
 If $d(A_1, A_2)> 2\delta$, let $p\in A_1$, $q\in A_2$ with
 $d(p,q)=d(A_1, A_2)$; in this case, let $B_1\subset A_1$
 be the closed segment with  midpoint $p$ and length $20\delta$
  and similarly for $B_2\subset A_2$. 
If $d(A_1, A_2)\le  2\delta$, 
let 
 $a=\min\{t\in \R: d(c(t), A_2)\le 2\delta\}$
 and $b=\max\{t\in \R: d(c(t), A_2)\le 2\delta\}$, and denote $y=c(a)$, 
$z=c(b)$; also let $y', z'\in A_2$ with $d(y, y')\le 2\delta$
  and $d(z,z')\le 2\delta$; in this case, let $B_1\subset A_1$ 
  be the  closed $20\delta$-neighborhood of $yz$  in   $A_1$
  and $B_2\subset A_2$ 
  be the  closed $20\delta$-neighborhood of $y'z'$  in   $A_2$.
We shall show that the two conditions in Lemma \ref{appping}
 are satisfied, and hence Proposition \ref{p2} holds.

\b{Le}\label{t1}
{Condition \e{(2)} of Lemma \ref{appping} is satisfied.}

\end{Le}

\b{proof}
We only write down the proof in the case 
$d(A_1, A_2)\le  2\delta$, $i=1$ and $n> 0$, the 
   other cases are similar or simpler.

Let $y\in p_1^{-1}(B_1)$.  
We need to show 
 $p_1(g_1^n(y))\notin B_1$.  Denote $z=p_1(y)$, $y'=g_1^n(y)$, 
 $z'=g_1^n(z)$  and $y''=p_1(y')$. Lemma \ref{pars}
  implies $d(z, s(z))\ge 160\delta$. Set $t=d(z, s(z))$.
 By Lemma \ref{order}
 and the fact that $d(p_1(s^i(z)), p_1(s^{i+1}(z)))\ge t-4\delta$ we obtain
 $d(p_1(z'), z)\ge 10k_1n(t-4\delta)$.

Let $z''=p_1(z')$. 
We claim  $d(y'', z'')\le 13 \delta$.
We first finish the proof assuming the claim. 
The claim implies $d(z, y'')\ge 10k_1(t-4\delta)n-13\delta$.
Note Lemma \ref{pars}
implies $t-4\delta\ge 2d(x, s(x))/5$. It follows that 
$d(z, y'')\ge 4nk_1d(x, s(x))-13\delta\ge 
4k_1d(x, s(x))-13\delta> 3k_1d(x, s(x))+40\delta\ge {\text{length}}(B_1)$.
  The last inequality follows from the definition of $B_1$ and Lemma \ref{twoax}.
  Since $z=p_1(y)\in B_1$, we have 
   $p_1(g_1^n(y))=y''\notin B_1$.

We now prove the claim. We may assume $d(y'', z'')>3 \delta$. 
  Fix geodesic segments
 $z'z''$, $y'z''$ and $y'y''$.
 Let $p\in y''z''$ with $d(p, y'')=3\delta$.
We use the $\delta$-thin condition. 
Consider the geodesic triangle $(y'y''z'')$. 
There is some $q\in y'y''\cup y'z''$ with
 $d(p,q)\le \delta$. 
 If $q\in y'y''$, then $d(q, y'')\ge d(p, y'')-d(p,q)\ge 2\delta$ and hence
  $d(q, p)< d(q, y'')$, contradicting the fact that $y''=p_1(y')$. 
 So $q\in y'z''$. Now consider the geodesic triangle $(y'z''z')$.
There is some $r\in z'z''\cup z'y'$ with $d(q,r)\le \delta$.  
If $r\in z'z''$, then $d(r, z'')\le d(z', z'')\le 2\delta$
  and hence 
 $d(y'',z'')=d(y'',p)+d(p, z'')\le 3\delta+ d(p,q)+d(q, r)+d(r, z'')\le 7\delta$. 
Suppose $r\in y'z'$.  Then $d(r, z')=d(r, g_1^n(A_1))\le d(r, p)+d(p,g_1^n(A_1) )
\le d(r, q)+d(q, p)+2\delta\le 4\delta$. 
 Therefore $d(z'', p)\le d(z'',z')+d(z', r)+d(r, p)\le 8\delta$ 
  and $d(z'', y'')\le 11\delta$.

\end{proof}

The following lemma can be proved by using approximating trees (\cite{GhD}).
%follows from the tree-like property of hyperbolic spaces???
One can also prove it directly from the $\delta$-thin condition. 
 Since it is more or less clear, 
  we omit  its  proof here.

\b{Le}\label{t2}
{Condition \e{(1)} of Lemma \ref{appping} is satisfied.}

\end{Le}

\qed

\subsection{Existence of   hyperbolic elements }\label{parab}

In this section we 
establish the following result, which guarantees the 
 existence of hyperbolic elements with short
  word length.

%shall show the existence of  a uniform hyperbolic element.

\b{Prop}\label{p3}
{Let $G$ be an infinite relatively hyperbolic group
  and $X$ a $\delta$-hyperbolic geodesic metric space  that $G$ acts upon as in 
    Definition \ref{rel}. Then 
there exists a positive integer $n_0$ with the following property:
   for any finite generating set $S$ of $G$, 
$S(n_0)$ contains a hyperbolic element.}

\end{Prop}

For a finite set of isometries $F$ of a metric space $X$, and $x\in X$,
 we let $\lambda(x,F)=\max\{d(f(x), x): f\in F\}$.   
 We use the following result of M. Koubi (\cite{K}).

\b{Prop}\label{Ko}
{Let 
$X$ be a $\delta$-hyperbolic geodesic metric space,
 and $G$ a group of isometries of $X$ with finite generating set
 $S$.  If $\lambda(x,S)>100\delta$ for all $x\in X$, then 
 $G$ contains a hyperbolic element  $g$ such that $d_S(id, g)=1$ or $2$.}

\end{Prop}

 Proposition \ref{p3}
 follows from
 Proposition \ref{Ko} and the following result.

\b{Le}\label{infin}
{There exists a positive integer $n_0$ with the following property.
 For any finite generating set $S$ of $G$, the inequality 
$\lambda(x,S(n_0))>100\delta$ holds for all $x\in X$.}

\end{Le}

\b{proof}
Let $\cal{B}$ be a $200\delta$-separated invariant system of horoballs centered at the
 parabolic points. Recall that 
 the action of $G$ on $Y(\cal{B})$ is proper and cocompact.
 Let $K\subset Y(\cal{B})$ be a compact set such that $\bigcup_{g\in G}g(K)=Y(\cal{B})$.
Set $a=\text{diam}(K)$ and fix a point $p\in K$. 
   Let $A=\{g\in G: d(g(p), p)\le 2a+100\delta\}$. Notice that $A$ is a finite set.
Since $G$ is infinite, there is  some $g_0\in G$ with 
 $d(g_0(p), p)>100\delta+2a$.   %Since $A$ is a finite set, 
  Notice that if $A$ generates $G$, 
 then 
$$n_0:=\max\{d_S(id, g_0): S\subset A \;\;\text{and}\;\; S \;\;\text{generates}\;\; G\}$$
  is finite.

Now let $S$ be a finite generating set of $G$. 
Suppose there is some $x\in X$ with 
$\lambda(x,S)\le 100\delta$. 
 By the definition of $\lambda(x,S)$, 
 we have $d(s(x), x)\le 100\delta$ for all $s\in S$. 
Recall that  $\cal{B}$ is  $200\delta$-separated. 
If $x\in B$ for some $B\in \cal{B}$, then $s(B)=B$ for all $s\in S$ and consequently
 the center of $B$  is fixed by the entire group $G$, a  contradiction. 
 Hence $x\in Y(\cal{B})$.  
There is some $\ga\in G$ with $\ga(x)\in K$. 
  For $s\in S$,  we have   $d(p, \ga s \ga^{-1}(p))\le d(p, \ga(x))+d(\ga(x), \ga s (x))
+d(\ga s (x), \ga s \ga^{-1}(p))=d(p, \ga(x))+d(x,  s (x))
+d(\ga (x), p)\le 2a+100\delta$.  
It follows that $\ga S\ga^{-1}:=\{\ga s\ga^{-1}:  s\in S\}\subset A$. 
 Clearly  $\ga S\ga^{-1}$ generates $G$.  There is some integer $k$, 
 $1\le k\le n_0$ and  $s_i\in S\cup S^{-1}$ ($1\le i\le  k$)  such that 
$g_0=(\ga s_1\ga^{-1})\cdots (\ga s_k\ga^{-1})=\ga (s_1\cdots s_k)\ga^{-1}$.
 %with $s_i\in S\cup S^{-1}$.
 Now $d(\ga^{-1}g_0\ga(x), x)=d(g_0\ga(x), \ga(x))\ge d(g_0(p), p)-d(g_0(p), g_0\ga(x))
-d(\ga(x), p)=d(g_0(p), p)-d(p, \ga(x)) 
-d(\ga(x), p)>100\delta+2a-a-a=100\delta$.
Notice $\ga^{-1}g_0\ga=s_1\cdots s_k\in S(n_0)$. 
 It follows that $\lambda(x, S(n_0))>100\delta$.

\end{proof}

 \addcontentsline{toc}{subsection}{References}
%\noindent 
%Xiangdong Xie\\
%Department of Mathematics\\
%Washington University \\
%St.Louis, MO 63130\\
%xxie@math.wustl.edu\\

\end{document}